\documentstyle[12pt]{amsart}
\newtheorem{theo}{THEOREM}[section]
\newtheorem{lemma}[theo]{Lemma}
\newtheorem{cor}[theo]{Corollary}
\newtheorem{prop}[theo]{Proposition}
\theoremstyle{definition}
\newtheorem{dfntn}[theo]{Definition}

\theoremstyle{remark}

\newtheorem{rem}[theo]{Remark}
\newtheorem{case}{Case} 

\newcommand{\brref}[1]{(\ref{#1})}
\newcommand{\tensor}{\otimes}
\newcommand{\bali}{\author[A. Alzati]{Alberto Alzati$^*$ }
\author[G. Besana]{ Gian Mario Besana$^{**}$}}
\newcommand{\Proj}[1]{\Bbb{ P}(#1)}
\newcommand{\Projcal}[1]{\Bbb{ P}({\cal #1})}

\newcommand{\restrict}[2]{{#1}_{\mid _{#2}}}

\newcommand{\calo}{{\cal O}}
\newcommand{\oof}[2]{{\cal O}_{#1}({#2})}
\newcommand{\oofp}[2]{{\cal O}_{\Bbb{ P}^{#1}}({#2})}
\newcommand{\iof}[2]{{\cal I}_{#1}({#2})}
\newcommand{\iofo}[1]{{\cal I}_{#1}}
\newcommand{\xel}{(X, L)}
\newcommand{\sel}{(S, L)}
\newcommand{\Pin}[1]{\Bbb{P}^{#1}}
\newcommand{\scroll}[1]{(\Proj{#1},\taut{#1})}
\newcommand{\taut}[1]{{\cal O}_{\Bbb{P}(#1)}(1)}
\newcommand{\tautof}[2]{{\cal O}_{\Bbb{P}(#1)}(#2)}
\newcommand{\tautcal}[1]{{\cal O}_{\Bbb{P}({\cal#1})}(1)}

\newcommand{\clc}{(C,\restrict{L}{C})}
\newcommand{\sls}{(S,\restrict{L}{S})}
\newcommand{\num}{\equiv}
\newcommand{\xt}{(X,\cal{T})}
\title[On the $k$-regularity]{On the $k$-regularity of
some projective manifolds }
\bali
\thanks{ $^*$ Dipartimento di Matematica - Via Saldini 50 - Universit\`a
degli Studi
di Milano - 20133 Milano - ITALY -  alzati@@vmimat.mat.unimi.it -
Supported by the
GNSAGAof the Italian C.N.R. ;\\
 $^{**}$ Department of Mathematics - Eastern Michigan University -
Ypsilanti MI -
U.S.A. 48197 - gbesana@@emunix.emich.edu.}
\begin{document}
\maketitle
\begin{center}
In memory of  F. Serrano
\end{center}

\begin{abstract} The conjecture on the $(degree - codimension +
1)$- regularity of projective varieties is  proved for smooth linearly
normal
polarized varieties $\xel$ with $L$ very ample,  for  low values of 
$\Delta \xel = degree  -  codimension -1.$ Results concerning the
projective
normality of some classes of special varieties including scrolls over
curves of genus
$2$ and quadric fibrations over elliptic curves, are proved.
\end{abstract}
%
%
\section{Introduction}

A complex projective variety $X\subset \Pin{N}$ is $k$-regular in the
sense of
Castelnuovo-Mumford if
$h^i(\iof{X}{k-i}) = 0$ for all
$i\ge 1$ where $\iofo{X}$ is the ideal sheaf of $X.$ If
$X$ is $k$-regular then the minimal generators of its homogeneous ideal
have
degree less than or equal to $k.$ 
A long standing conjecture, known to us as the Eisenbud Goto conjecture,
states
that an $n$-dimensional  variety $X \subset \Pin{N}$ of degree $\deg X =
d$ is
$(d-(N-n)+1)$-regular. Gruson Lazarsfeld and Peskine \cite{GLP}
established the
conjecture for curves, Lazarsfeld \cite{Laz1} for smooth surfaces and
Ran
\cite{Ran1} for threefolds with high enough codimension. A nice
historical account
of the conjecture and further results can be found in \cite{Kw}.
 In section \ref{deltacongettura} the conjecture is proved for all
smooth
varieties $X$  embedded by the complete linear system associated with a
very
ample line bundle $L$ such that $\Delta \xel \le 5$ where $\Delta \xel =
\dim{X} +
\deg{X} -h^0(L).$
 Notice also that in  recent times computer algebra systems like
 Macaulay have made possible the
explicit construction and study of examples of algebraic varieties
starting from
minimal generators of the homogeneous ideal of the variety. A priori
information
on the $k$-regularity of a variety is therefore useful for these
constructions.

Strictly related to the notion of $k$-regularity is the notion of
$k$-normality of a
projective variety. A variety $X \subset \Pin{N}$ is $k$-normal if
hypersurfaces of
degree $k$ cut  a complete linear system on $X$ or, equivalently, 
if $h^1(\iof{X}{k})=0.$ If $X$ is $k$-regular it is clearly
$(k-1)$-normal. $X$ is said to
be  projectively normal if it is $k$-normal for all $k\ge 1.$

As a by-product of the proof of the above result  the projective
normality of a class of surfaces of degree nine in $\Pin{5}$ which was
left as an open
question in \cite{gisa} is established in Lemma \ref{blowupF1}. The non
existence of a
class of scrolls of degree $10$, left as an open problem in
\cite{fa-li10}, is also
established in Remark \ref{scr10}.

In section \ref{g2scrolls} we deal with the projective normality of
scrolls $X
=\Proj{E}$ over a curve of genus $2$ embedded by the complete linear
system
associated with the tautological line bundle $\taut{E},$ assumed to be
very ample.
Two-dimensional such scrolls are shown to be always projectively normal
except
for a class $S$ of non
$2$-normal surfaces of degree eight in
$\Pin{5}$ studied in detail in \cite{alibaba}. Three-dimensional scrolls
$X=\Proj{E}$
of degree
$\deg X \ge 13$ are then shown to be projectively normal if and only if
$E$ does not
admit a quotient
$E
\to
\cal{E}
\to 0$ where $P(\cal{E})$ belongs to the class $S$  of non quadratically
normal
surfaces mentioned above.

In section \ref{ellipticpkbundles}, building on the work of Homma
\cite{Ho1},\cite{Ho2} and Purnaprajna and Gallego \cite{pu-ga}, criteria
for
the projective normality of three-dimensional quadric bundles over
elliptic curves
are given, improving some results contained in \cite{bu}.
%
%
\section{General Results and Preliminaries}
\label{prelimsec}
\subsection{Notation}
\label{notation}
The notation used in this  work is mostly standard from Algebraic
Geometry. Good references are \cite{H} and \cite{gh}.
The ground field is always the field $\Bbb{ C}$ of complex
numbers. Unless otherwise stated all varieties are supposed to be
projective.
$\Bbb{P}^{N}$ denotes the N-dimensional complex projective space.
 Given a projective n-dimensional variety $X$, ${\cal O}_X$
denotes its structure sheaf  and  $Pic(X)$  denotes the group of line
bundles over $X.$ Line bundles, vector bundles and Cartier
divisors are denoted by capital letters as $L, M,\cal{M} \dots.$ Locally
free sheaves of
rank one, line bundles and Cartier divisors are used interchangeably as
customary. 

Let $L, M \in Pic(X)$, let $E$ be a vector bundle of rank  $r$ on $X$,
let
$\cal{F}$ be a coherent sheaf on $X$ and let $Y\subset X$ be a
subvariety of $X.$
Then the following notation is used:
\begin{enumerate}
\item[ ] $LM$ the intersection of divisors $L$ and $M$
\item[ ] $L^{n}$ the degree of $L,$
\item[ ] $|L|$ the complete linear system of effective divisors
associated with $L$,
\item[ ]$L_Y$ or $\restrict{L}{Y}$ the restriction of $L$ to  $Y,$
\item[ ] $L \sim M$  linear equivalence of divisors
\item[ ] $L \equiv M$ numerical equivalence of divisors
\item[ ] Num$(X)$ the group of line bundles on $X$ modulo  numerical
equivalence
\item[ ]  $\Bbb{P}(E)$ the projectivized bundle of $E,$ see \cite{H}
\item[ ] $H^i(X, \cal{F})$ the $i^{th}$ cohomology vector space with
coefficients in ${\cal
F},$ 
\item[ ] $h^i(X,\cal{F})$ the dimension of $H^i(X, \cal{F}),$ here and
immediately
above $X$ is sometimes omitted when no confusion arises.
 \end{enumerate}

If $C$ denotes a smooth projective curve of genus $ g$, and $E$ a vector
bundle
over $C$ of deg $E= c_1(E)= d$ and rk$E=r$, we need the following
standard
definitions:
\begin{enumerate}
\item[ ] $E$ is $\it normalized$ if $h^0(E)\ne 0$ and $h^0(E \otimes
\cal{L})=0$
for any invertible sheaf $\cal{L}$ over $C$ with deg$\cal{L}<0$.
\item[ ] $E$ has slope $\mu(E) = \frac{d}{r}$.
\item[ ] $E$ is $\it semistable$ if and only if for every proper
subbundle $S$,
$\mu(S) \leq \mu(E)$. It is $\it stable$ if and only if the inequality
is strict.
\item[ ] The Harder-Narasimhan filtration of $E$ is the
unique filtration:
$$0=E_0\subset E_1\subset ....\subset E_s=E$$
 such that
$\frac{E_i}{E_{i-1}}$ is semistable for all $i$, and $\mu_i(E)=\mu
(\frac{E_i}{E_{i-1}})$ is
a strictly  decreasing function of $\it i$.
\end{enumerate}
A few definitions from  \cite{bu} needed in the sequel are recalled.

Let  $0=E_0 \subset E_1 \subset ....\subset E_s=E$ be the
Harder-Narasimhan
filtration of a vector bundle $E$ over $C$. Then
\begin{enumerate}
\item[]$\mu^-(E)=\mu_s(E)=\mu (\frac{E_s}{E_{s-1}})$
\item[]$\mu^+(E)=\mu_1(E)=\mu (E_1)$
\item[]or alternatively
\item[]$\mu^+(E)= \text{ max }\{\mu(S) |0 \to S \to E \}$
\item[]$\mu^-(E)= \text{ min }\{\mu(Q) |E \to Q \to 0 \}$.
\end{enumerate}
It is  also $\mu^+(E) \geq \mu(E) \geq \mu^-(E)$ with equality if and
only if $E$ is
semistable.
In particular if $C$ is an elliptic curve, an indecomposable vector
bundle $E$ on $C$ is
semistable and hence  $\mu(E) = \mu^-(E)=\mu^+(E)$.

 The following definitions are standard in the theory of polarized
varieties. A good
reference is \cite{fu}. A {\em polarized variety } is a pair $\xel$
where $X$
is a smooth projective n-dimensional variety and $L$ is an
ample line bundle on $X$. Its {\em  sectional genus}, denoted
$g\xel$, is defined by $2g\xel - 2 = (K_X +
(n-1)L)  L^{ n-1}$. Given any $n$-dimensional polarized
variety $\xel$ its $\Delta${\em - genus }  is defined by $\Delta \xel =
dim (X) + L^{ n} - h^0\xel.$
A polarized variety $(X, L)$ has a {\em ladder} if  there
exists a sequence of reduced and irreducible subvarieties $X = X_n
\supset X_{n-1} \dots \supset X_1$ of $X$ where $X_j\in|L_{j+1}| =
|\restrict{L}{X_{j+1}}|.$ Each $(X_j, L_j)$ is called a   {\em rung} of
the ladder. If $L$ is generated by global sections $(X, L)$ has a
ladder. A rung $(X_j,
L_j)$ is {\em regular}  if
$H^0(X_{j+1},\restrict{L}{X_{j+1}})
\to H^0(X_j,\restrict{L}{X_j})$ is onto. The ladder is regular if all
the
rungs are
regular. If the ladder is regular $\Delta (X_j,L_j) = \Delta(X,L)$ for
all $1\le j \le n.$
A variety $X \subset \Pin{N}$ is {\it k-normal} for some $k \in \Bbb{Z}$
if
$H^0(\Pin{N},
\oofp{N}{k})
\to H^0(\oof{X}{k})$ is onto. Equivalently, if $\iofo{X}$ is the ideal
sheaf of $X,$ $X$ is
$k$-normal if
$h^1(\iof{X}{k}) = 0.$
$X$ is {\it projectively normal} if it is $k$-normal for all $k \ge 1.$
A polarized pair
$\xel$ with $L$ very ample is called $k$-normal or projectively normal
if $X$ is
$k$-normal or p.n. in the embedding given by $|L|.$ A polarized variety 
$(X, L)$ with
$L$ very ample is always $1$-normal (linearly normal). 

 A line bundle $L$ on
$X$ is {\em normally (or simply)  generated} if the graded algebra
$G(X, L) =  \bigoplus_{t\ge 0}H^0(X,tL)$ is generated by $H^0(X,L).$ 
$L$
is very ample and normally generated if and only if   $(X, L)$  is p.n.

A variety $X\subset \Pin{N}$ is $k$-regular, in the sense of
Castelnuovo-Mumford,  if
for all
$i\ge 1$ it is $h^i(\iof{X}{k-i})=0.$  A polarized pair $\xel$ with $L$
very ample is
$k$-regular if 
$X$ is $k$-regular in the embedding given by $|L|.$
If $X$ is $k$-regular then it is $(k+1)$-regular.
\medskip
\subsection{General Results}
\medskip
Let $C$ be a smooth  projective curve of genus $g$, $E$ a vector
bundle of rank $n$, with $n \ge 2$, over $C$ and $\pi : X =\Bbb{P}(E)
 \to C $ the projectivized bundle
associated to $E$ with the natural projection $\pi$.
 Denote
with $\cal{T} = \taut{E} $ the tautological sheaf and with  $\frak{F}_P=
\pi^*\cal {O}_{C}(P) $  the line bundle associated with the fiber over
$P\in C.$ Let
$T$ and
$F$ denote the numerical classes respectively of $\cal T$ and
$\frak{F}_P$.
In this work we refer to  a  polarized variety $\xt$ as a  {\it scroll}
over a curve $C$ if
there is a vector bundle
$E$ over $C$ such that $\xt= \scroll{E}$  and $\cal{T}$ is very ample.
\begin{rem}
\label{leray}
Let $\pi : \scroll{E} \to C$ be a $n$-dimensional projectivized bundle
over a curve
$C.$ From Leray's Spectral sequence and standard facts about higher
direct image sheaves (see for example
\cite{H} pg. 253) it follows that 
\begin{gather}
H^1(\tautof{E}{t})=H^1(C, S^tE) \text{ for } t\ge 0 \notag\\ 
H^i(\tautof{E}{t}) = 0 \text {
for } i \ge 2  \text { and } t>-n. \notag
\end{gather}
\end{rem}
 
Let $D\sim a\cal{T} + \pi^*B$, with  $ a\in \Bbb{Z}$, $B\in Pic(C)$
and $\deg B = b$, then $  D \equiv aT+bF .$  Moreover 
$\pi_*(\oof{\Proj{E}}{D}) =
S^{a}(E)
\otimes
\cal {O}_{C}(B) $ and hence $\mu^-(\pi_*(\oof{\Proj{E}}{D})=a\mu^-(E)
+b$ (see
\cite{bu}).

 Regarding  the ampleness, the global generation, and
the normal generation  of $D$, a few known criteria useful in the sequel
are listed
here:
 \begin{theo}[Miyaoka \cite{Miyao3}]
\label{miyaoteo}
 Let $E$ be a vector bundle over a smooth
projective curve $C$ of genus $g$, and $X =\Bbb{P}(E)$ . If $D\equiv
aT+bF$ is a line
bundle over $X$, then $D$ is ample if and only if $a>0$ and $b+a
\mu^-(E) >0$.
\end {theo}

\begin{lemma}
\label{buongg}
(see e.g. \cite{bu}, Lemma 1.12)
Let $E$ be a vector bundle over $C$ of genus $g$.
\begin{itemize}
\item[i)]if $\mu^-(E) > 2g-2$ then $h^1(C, E)=0$
\item[ii)]if $\mu^-(E) > 2g-1$ then $E$ is generated by global sections.
\end{itemize}
\end{lemma}

\begin{lemma}[Butler \cite{bu} Theorem 5.1A]
\label{criteriodelbutler}
 Let $E$ be a vector bundle on a
smooth projective curve of genus $g$ and let $D\num aT + bF$ be a
divisor on $X
=\Bbb{P}(E).$ If
\begin{equation}
\label{condizionedelbutler}
b+a\mu^-(E) > 2g.
\end{equation} 
 then $D$ is normally generated.
\end{lemma}

%
%
A few basic facts on the {\it Clifford Index} of a curve are recalled.
Good references
are  \cite{mart} and  \cite{GL}. Let $C$ be a projective curve and
$L$ be any line bundle on $C$. The Clifford index of $L$  is defined as
follows:
$$cl(L)=deg(L) -2(h^0(L)-1).$$
The Clifford index of the curve is  $cl(C)=\text {min}\{cl(L) |
h^0(L)\geq 2 \text{ and }h^1(L)\geq 2 \}$.
 For a general curve $C$
it is $cl(C)=\left [\frac{g-1}{2}\right ]$
 and in any case $cl(C)\leq\left [\frac{g-1}{2}\right ]$. By
Clifford's  theorem a  special line bundle $L$ on $C$ has $cl(L)\geq
0$ and
the equality holds if and only if $C$ is hyperelliptic and $L$ is a
multiple of the unique
$g^1_2$.\\ If
$cl(C)=1$ then $C$ is  either a plane quintic curve or a trigonal curve.
\begin{theo}[\cite{GL}]
\label{glcliff}
Let L be a very ample line bundle on a smooth irreducible complex
projective curve
$C.$ If 
$$\deg(L)\geq 2g+1-2h^1(L)-cl(C)$$
then $(C, L)$ is projectively normal.
\end{theo}

\section{The Eisenbud Goto conjecture for low values of $\Delta.$}
\label{deltacongettura}
Let $X \subset \Pin{N}$ be an $n$ dimensional projective variety of
degree $d$. A long
standing conjecture, known to us as the {\it Eisenbud Goto} conjecture,
states that
$X$ should be $(d-(N-n) +1)$-regular, i.e.
$(degree-codimension+1)$-regular.

Many authors worked on the conjecture for  low values of the dimension
and
codimension of X. A nice historic account is found in \cite{Kw}. Some of
their results
are collected in the following Theorem.
\begin{theo}[\cite{GLP}, \cite{Laz1}]
\label{lazran}
If $X\subset \Pin{N}$ is any smooth  curve or  any smooth surface then
$X$ is 
$(d-c+1)$-regular where $d =\deg{(X)}$ and $c = \text{codimension} (X)$.
\end{theo}

 In this section we would like to offer a proof of the conjecture for
linearly normal smooth varieties with low $\Delta$-genus.   Let
$\xel$ be a polarized variety with
$L$ very ample. The above conjecture can be restated for the embedding
given by
$|L|$ in terms of
$\Delta$-genus as follows:

{\bf Conjecture } {\it Let $\xel$ be a polarized variety with $L$ very
ample. Then $\xel$
is
$(\Delta + 2)$-regular.}
%
%
\begin{rem}
\label{ipers}
It is  straightforward  to check  that hypersurfaces of degree $d$ are
always
$d$-regular and not $(d-1)$-regular. This shows that the conjecture is
indeed
sharp. On the other hand  there are varieties $X\subset \Pin{N}$ which
are
$k$-regular for $k < d-c+1.$ This motivates Definition \ref{extremal}. 
\end{rem}

\begin{rem}
It is  a classical adjunction theoretic results that given $\xel$ with
$L$ very ample, $K_X + tL$ is globally generated, and in particular
$h^0(K_X + tL) \ne 0,$ for
$ t
\ge n$ unless
$t=n$ and $\xel = (\Pin{n}, \oofp{n}{1}).$ This fact, Remark \ref{ipers}
and the sequence $0 \to \iofo{X}\to \calo_{\Pin{N}} \to \calo_X \to 0$
suitably twisted  show that no linearly normal non degenerate
$n$-dimensional  variety $X \subset \Pin{N}$ can be $k$-regular for
$k \le 1.$ Therefore in what follows we will always assume $k \ge 2$
when dealing with $k$-regularity.
\end{rem}
\begin{dfntn} 
\label{extremal}  Let $X \subset
\Pin{N}$ be  a $n$-dimensional   variety   of degree
$d.$Let
$$r(X)= Min \{k \in \Bbb{Z} | X \text{ is } k-\text{regular}\}.$$  A
variety   $X$ is
{\bf extremal} if
$r (X) = d-(N-n)-1.$ 
A polarized variety $\xel$ with $L$ very ample is {\bf extremal} if it
is extremal in
the embedding given by $|L|,$ i.e. if $r(X,L) = Min
\{k \in \Bbb{Z} | \xel \text{ is } k-\text{regular}\} = \Delta + 2.$
\end{dfntn}
In what follows we will prove  the above conjecture for all
linearly normal manifolds with 
$\Delta
\le 5$ obtaining along the way the value of $r \xel$ for most  of  the
same
manifolds.

%
%
\begin{lemma}
\label{hyperplanesec}
Let $X \subset \Pin{N}$ be a smooth n-dimensional variety and let   $Y
\subset
\Pin{N-1}$  be a generic hyperplane section. 
\begin{itemize}
\item[i)] If $X$ is $k$-regular then $Y$ is $k$-regular
\item[ii)] If $Y$ is $k$-regular and $X$ is $(k-1)$-normal then $X$ is
$k$-regular.
\item[iii)] If $X$ is $(r(Y)-1)$-normal then $r(X)=r(Y)$
\end{itemize}
\end{lemma}
\begin{pf}
The exact sequence
\begin{equation}
\label{ideali}
0 \to \iof{X}{k-i} \to \iof{X}{k-i+1} \to\iof{Y}{k-i+1} \to 0.
\end{equation}
immediately gives i). To see ii) consider again sequence \brref{ideali}.
The $k$-regularity of $Y$ gives
$h^{i-1}(\iof{Y}{k-i+1})=0$ for all
$i\ge 2.$ Since $k$ regularity implies $k+1$-regularity it is
$h^i(\iof{Y}{k-i+1})=0$ for all $i \ge 1.$ Therefore 
$h^i(\iof{X}{k-i})=h^i(\iof{X}{k-i+1})$  for all $i\ge 2$ from
\brref{ideali} and
iteratively
$h^i(\iof{X}{k-i})=h^i(\iof{X}{k-i+t}$ for all $i\ge 2$ and for all
$t\ge 1.$ Letting $t$
grow, Serre's vanishing theorem gives $h^i(\iof{X}{k-i+t}=0$ for all
$i\ge 2$ and all
$t\ge 1$ and thus 
 $h^i(\iof{X}{k-i})=0$ for all $i\ge 2.$ Because $X$ is assumed
$(k-1)$-normal it is
$h^1(\iof{X}{k-1})=0$ which concludes the proof of $ii).$
Now $iii)$ follows immediately from $i)$ and $ii).$
\end{pf}
 
\begin{lemma}
\label{hyperplanesecpol}
Let  $\xel$ be a polarized variety with $L$ very ample. Let $Y \in |L|$
be a generic
element and assume 
$H^0(X, L)\to H^0(Y,  \restrict{L}{Y})$ is onto. Then $i), ii),iii)$ as
in Lemma
\ref{hyperplanesec} hold if we replace $X$ by $\xel$ and $Y$ by $(Y,
\restrict{L}{Y}).$
\end{lemma}
\begin{pf}
Let $h^0(L)=N+1.$ The surjectivity condition on the restriction map
between global
sections of
$L$ and
$\restrict{L}{Y}$ guarantees that $|\restrict{L}{Y}|$ embedds $Y$ as a
linearly normal
manifold in $\Pin{N-1}$ , therefore the same proof as in Lemma
\ref{hyperplanesec} applies.
\end{pf}

\begin{rem} 
\label{liftingpn}
Let  $\xel$ be a polarized variety with $L$ very ample. Let $Y \in |L|$
be a generic element and assume 
$H^0(X, L)\to H^0(Y,  \restrict{L}{Y})$ is onto. Then \cite{fu}
Corollary 2.5 shows
that if $(Y, \restrict{L}{Y})$ is projectively normal, so is $\xel.$
Therefore when the
ladder is regular and $Y$ is p.n. Lemma \ref{hyperplanesecpol} gives $r
\xel = r (Y,
\restrict{L}{Y}).$
\end{rem}

\begin{lemma}
\label{rofpncurves}
Let $(C, L)$ be a projectively normal curve with $g\ge 1.$ 

Then $r (C, L) = Min \{ t \ge 3 |
h^1((t-2) L) =0\}.$
\end{lemma}
\begin{pf}
Let $h^0(L)=N+1$ so that $C\subset \Pin{N}.$ It is 
$h^1(\iof{C}{k-1})=0$ for all
$k\ge 2$ because of the projective normality assumption. The sequence
$$0 \to \iof{C}{k-i} \to \oofp{N}{k-i} \to (k-i) L \to 0$$
easily gives $h^i(\iof{C}{k-i})=0$ for all $i \ge 3$ and $k \ge 2.$ 

The same sequence
gives $h^2(\iof{C}{k-2})=h^1((k-2) L)$ and since 
$h^1(\calo_C) = g \ge 1$ it is $r (C, L) = Min \{ t \ge 3 | h^1((t-2) L)
=0\}.$
\end{pf}

In order to apply the above lemmata in one occasion
the projective normality of a particular class of
surfaces of degree nine needs to be established. The
following Lemma also improves \cite{gisa}. Here $\Bbb{F}_1$ denotes the
Hirzebruch
rational ruled surface of invariant $e=1,$  $\pi : Bl_{t}S \to S$
denotes the blow up of
a surface $S$ at $t$ points, $E_i$ are the exceptional divisors of the
blow up, 
$\frak{C}_0=\pi^*(C_0)$ denotes the pull back of the line bundle
associated with the
fundamental section of
$\Bbb{F}_1$ and
$\frak{f}=\pi^*(f)$ the pull back of the one associated with  any fibre
$f$ of the
natural projection
$p:
\Bbb{F}_1
\to \Pin{1}.$
\begin{lemma}
\label{blowupF1}
Let $\sel=(Bl_{12} {\Bbb F}_1, 3{\frak C_0} + 5{\frak f} - \sum_i E_i)$.
Then $\sel$ is projectively normal.
\end{lemma}
\begin{pf} 
The projective normality of linearly normal degree nine surfaces
was studied in \cite{gisa}. Let $\sel$ be a surface of degree $9$ and
sectional genus
$5,$ embedded in $\Pin{5}.$
 The surface under consideration was
established to be projectively normal unless its generic curve
section $C$ is  trigonal  and $\restrict{L}{C}=K_C-M+D$ where $M$ is a
divisor in the
$g^1_3$ and $D$ is a divisor of degree $4$ giving a foursecant line for
$C.$
Therefore if $S$ were not p.n. it would admit an infinite number of
$k\ge 4$-secant lines. On the other hand a  careful study of the
embedding shows
that
$S$ contains only  a finite number of lines and that the only lines 
with
 self intersection $\ge -1$ are the $12$ exceptional divisors $E_i.$
Thus the formulas
contained in
\cite{LeBz1} can be used.  A straightforward calculation using
\cite{LeBz1} shows that $S$ cannot have a infinite number of $k\ge
4$-secants,
contradiction.
\end{pf}
%
%

\begin{theo}
\label{alafujita}
Let $\xel$ be a $n$ dimensional polarized pair, $n\ge 2,$ with a ladder. 
Assume $g=g \xel \ge \Delta \xel = \Delta$ and $d=L^n \ge 2 \Delta + 1.$
Then :
\begin{itemize}
\item[i)] The
curve section  $\clc$ is $k$-regular if and only if $\xel$ is
$k$-regular and $r \xel= r
(C,\restrict{L}{C}).$
\item[ii)] Either $\Delta=0,1$ and $\xel$ is extremal or $\Delta\ge 2$
and $r \xel=3.$
\end{itemize}
\end{theo}
\begin{pf}
From \cite{fu} Theorem
(3.5) and from the fact that a normally genereated ample line bundle is
automatically
very ample it follows that
$L$ is very ample, $g=\Delta$, the ladder is regular and every rung of
the ladder  is
projectively normal. Therefore Lemma \ref{hyperplanesecpol} immediately
gives
$i).$ 

Then  $\xel$ is extremal if and only if the curve section $\clc$ is
such. Extremal linearly normal curves were classified in \cite{GLP} and
they are
either rational or elliptic normal curves. Therefore
$\xel$ is extremal if and only if $\Delta=g=0,1.$
Now assume $\Delta \ge 2$ and thus $\xel$ not extremal. The curve
section $\clc$ is
embedded in $\Pin{M}$where $M= d - \Delta.$  Since
$h^M(\oofp{M}{k-M})=h^0(\oofp{M}{-1-k})=0$ for all $k \ge 0$ the
sequence $0 \to
\iof{C}{k-i} \to \oofp{M}{k-i} \to \oof{C}{k-i} \to 0$ shows that
$h^i(\iof{C}{k-i}) =
h^{i-1}(\oof{C}{k-i})$ for all $i \ge 2$ and all $k\ge 0.$
Therefore, because $h^1(\calo_C) = g \ge 2,$ it must be $r \xel \ge 3.$
If $i\ge 3$ then clearly $h^{i-1}(\oof{C}{3-i}) =0$ and thus
$h^i(\iof{C}{3-i})=0.$
It is also $h^2(\iof{C}{1}) = h^1(\oof{C}{1}) = h^1(\restrict{L}{C})=0$
because
$g=\Delta$ and $d\ge 2 \Delta+1 > 2g-2.$
Since every rung of the ladder is projectively normal, in
particular $h^1(\iof{C}{2}) =0$ and thus $\clc$ is  $3$-regular. We can
conclude
that $r
\xel=r
\clc = 3.$
\end{pf}
\begin{prop}
\label{elscr}
Let $\xt = (\Proj{E}, \taut{E}) $ be a scroll over an elliptic curve.
Then
$r \xt = 3.$
\end{prop}
\begin{pf}
Because $h^2(\iofo{X})=h^1(\calo_X)=1$ it is $r \xel \ge 3.$
We need to show that $h^i(\iof{X}{3-i})=0$ for all $i\ge 1.$ Notice that
$|\cal{T}|$ embeds $X$ into $\Pin{N}$ as a variety of degree $d$ where
$N=d-1.$  Let  $i=1.$ It is known, cf. \cite{bu} and \cite{alibaba},
that
elliptic scrolls are projectively normal, so
$ h^1(\iof{X}{2}) = 0.$
Let $i=2.$ From Remark \ref{leray} it is $h^2(\iof{X}{1}) =
h^1(\oof{X}{1}) = h^1(C, E).$ Because $E$ is very ample it is
$\mu^-(E)
>0$ which, by Lemma \ref{buongg} implies  $h^1(C, E) = 0.$

  For
 $i=N$ it is
$h^N( \oofp{N}{3-N}) = h^0(\oofp{N}{-4})=0.$ 
Therefore it follows that
$h^i(\iof{X}{3-i}) = h^{i-1}(\oof{X}{3-i})$ for all $i\ge 3.$
Remark \ref{leray} gives
$h^{i-1}(\oof{X}{3-i})  = 0$ for  $3\le  i \le n+1$ while clearly
$h^{i-1}(\oof{X}{3-i})=0$ for $i>n + 1$ since $n=\dim X.$ 

Therefore 
$h^i(\iof{X}{3-i}) = 0$ for all $i\ge 3.$
\end{pf}
\begin{lemma}
\label{diseqonscr}
Let $(\Proj{E},\taut{E})$ 
 be a n-dimensional scroll over a curve of genus $g \ge 2.$  If
$\deg (E) >2g-2$ then
$\Delta
\ge 2n+g-3.$
\end{lemma}
\begin{pf}
Because $d=\deg(det E)=\deg(E) > 2g-2$, it is  $h^0(det E) = 1 + d - g$
by Riemann
Roch. Combining this with the inequality $h^0(det E) \ge h^0(E) + r-2$
found in
\cite{Io-To}, it follows that $h^0(E) \le d - n + 3-g$ and therefore
$\Delta \ge 2n + g
-3.$
\end{pf}
\begin{rem}
\label{scr10}
Notice that the above Lemma \ref{diseqonscr} rules out the existence of
scrolls of
degree
$10$ over a curve of genus $g=3$ left as an open possibility in
\cite{fa-li10}.
\end{rem}
%
We can now prove the main theorems of this section. For $\Delta \le 3$
we
establish the conjecture and give the value of $r \xel$ for all pairs.
For $\Delta=4,5$
we establish the conjecture and collect in a remark the known values of
$r \xel.$

\begin{theo}
\label{Delta+2thm}
Let $\xel$ be a n-dimensional  polarized pair with $X$ smooth, $L$ very
ample and
$\Delta
\le 5.$ Then $\xel$ is $\Delta + 2$-regular.
\end{theo}
\begin{pf}
Because of Theorem \ref{lazran} and Remark \ref{ipers}, the blanket
hypothesis
$n\ge 3$ and $codim X \ge 2$ will be in place throughout this proof.
\begin{case}
$\Delta \le 1$
\end{case}
If $\Delta = 0$ then $\xel$ is extremal by Theorem  \ref{alafujita}.
Assume
$\Delta=1$, because
$g=0$ implies
$\Delta=0,$ see \cite{fu} Prop. (3.4),  it is $g \ge1.$ Because $\xel$
is not a
hypersurface it is 
$d\ge 3$  and again Theorem
\ref{alafujita} gives $\xel$ extremal. 
\begin{case}
$\Delta=2$
\end{case}
If $g\le 1$ then
$\xel$ must be a two dimensional elliptic scroll, see \cite{fu} Theorem
(10.2).
Proposition
\ref{elscr} gives
$r=3.$ Let  $g\ge 2.$  Because $X$ is not a hypersurface  it is 
$h^0(L) \ge n+3.$ This implies $\Delta\le d-3$ i.e. $d\ge 5$ and then 
Theorem
\ref{alafujita} gives $r=3.$  
\begin{case}
$\Delta=3$
\end{case}
From
\cite{BEL2} it follows that  complete intersections of type
$(2,3)$ have $r=4.$  Following \cite{Io1}
Theorem 4.8 and section 7, it follows from Theorem \ref{alafujita} and
Proposition
\ref{elscr} that  the only varieties  left to investigate are Bordiga 
threefolds scrolls in $\Pin{5}.$ They  have  the following resolution
with
 $N=5.$ 
\begin{equation}
\label{bordigares}
0 \to \oofp{N}{-4}^{\oplus 3} \to \oofp{N}{-3}^{\oplus 4} \to \iofo{X}
\to 0.
\end{equation}
 Equalities $h^{N-1}(\iof{X}{3-N}) =
h^0(\calo_{\Pin{N}})=1$  and
$h^i(\iof{X}{3-i})=0$ for all $i \ge 1$ are straightforward to see, 
therefore $r=3.$
\begin{case}
$\Delta=4$
\end{case}
Varieties with
$\Delta=4$ are classified. Let us follow the list of varieties
given  in  \cite{Io4} Theorem 3. Threefolds in $\Pin{5}$ with  $d=7$,
$g=5$ or 
$g=6$ have respective  resolutions as in \brref{D4g5res} and
\brref{D4g6res}
with $N=5.$  
\begin{equation}
\label{D4g5res}
0 \to \oofp{N}{-5} \oplus \oofp{N}{-4} \to \oofp{N}{-3}^{\oplus 3} \to
\iofo{X}
\to 0.
\end{equation}
\begin{equation}
\label{D4g6res}
0 \to \oofp{N}{-5}^{ \oplus 2} \to \oofp{N}{-2} \oplus
\oofp{N}{-4}^{\oplus 2} \to
\iofo{X}
\to 0.
\end{equation}
The resolutions
\brref{D4g5res}  and \brref{D4g6res} quickly show that
$r=4.$
Complete intersections of type
$(2,2,2)$ have
$r=4$ by
\cite{BEL2}. 
Scrolls over  a genus $2$ curve must be two-dimensional while elliptic
scrolls
are taken care of by Proposition
\ref{elscr}.

Let now $q=0$ and $g=4.$ If $d\ge 9$  Theorem \ref{alafujita} gives
$r=3.$  On the
other hand since $\Delta=4$ and the codimension must be at least two, it
follows that
$d\ge 7.$ Let us now compare the varieties under consideration with the
lists of
manifolds of degree $7$ and $8$ given in \cite{Io1} and \cite{Io2}.

If $d=8$  then
$X\subset
\Pin{6}$ is a threefold scroll over the quadric surface. Since
$q=0$ the ladder is regular. Consider the curve section
$(C,\restrict{L}{C}).$ 
Such a  $(C, \restrict{L}{C})$ is known to be non hyperelliptic (see
\cite{Io2})  and thus Theorem \ref{glcliff}  gives $(C,
\restrict{L}{C})$ p.n. Since
$d =8 > 2g-2 = 6$ it is $h^1(\restrict{L}{C})=0$ and thus $r(C,
\restrict{L}{C})=3$ by lemma \ref{rofpncurves}. Because the ladder is
regular 
$r \xel=3$ by Remark
\ref{liftingpn}. 

 If  $d=7$ then $X \subset \Pin{5}$ is Palatini's scroll
over the cubic surface. A resolution for $\iofo{X}$ is found in
\cite{BSS3}:
$$0 \to \oofp{5}{-4}^{\oplus 4} \to \Omega^1(-2) \to \iofo{X} \to 0.$$
A simple cohomological calculation gives $r=4.$

\begin{case}
$\Delta=5$
\end{case}
Theorem \ref{alafujita} takes care of cases with $g\ge 5$ and $d\ge 11.$
Manifolds
with degree $d\le 10$ were classified by various authors and we will
examine them
later in the proof. Let us now assume $d \ge 11$ and $g\le 4.$  
Because $\Delta=5$ and elliptic scrolls are dealt with in Proposition
\ref{elscr}, it
must be
$g \ge 2.$ Varieties of low sectional genus were classified in 
\cite{Io1}. Let us follow the lists given there. If $g=2$ scrolls over a
curve are the only
manifolds to be considered. On the other hand such scrolls of genus $2$
have
$\Delta=2n$ (cf.
\cite{gisa}) so there are no manifolds  to examine. If $g=3$ scrolls
over curves are
ruled out by Lemma \ref{diseqonscr} and   scroll over $\Pin{2},$ having
$q=0,$ are
ruled out by \cite{Io1} Theorem 4.8 iv).  If $g=4$ scrolls over curves
are again ruled
out by Lemma \ref{diseqonscr}. Using  standard numerical relations  (see
for
example \cite{fa-li9} ( 0.14)) one sees that there are no hyperquadric
fibrations of
dimension $n\ge 3,$ $g=4,$ $\Delta=5$ over
$\Pin{1}$ or over an elliptic curve. Let now $\xel$ be a threefold which
is a scroll over
a surface
$(Y, \cal{L})$ with $q(Y)=0,$ $g \xel=4.$  Because 
$h^1(\calo_X)=q(Y)=0,$  recalling that a general hyperplane section of
$X$ is
birational to $Y$ and thus regular, the  ladder is regular and then
$\Delta
\xel =
\Delta(C,
\restrict{L}{C}) = 4$ by Riemann Roch.

Let us now consider the cases with $d\le 10$  by looking at the
classification
found in \cite{Io2},  \cite{fa-li9},  \cite{fa-li10}. The first non
trivial case occurs with
$d=8.$ $\xel$ is a threefold in $\Pin{5},$ admitting a fibration over
$\Pin{1}$ with
generic  fibers complete intersections of type $(2,2)$ in $\Pin{4}.$ A
resolution of the
ideal of this variety can be found in \cite{BSS3}. A standard
cohomological calculation
shows that  $r \xel = 4.$

Let now $d=9.$ From \cite{fa-li9} all varieties to be considered are
threefolds in
$\Pin{6}$ with 
$g=5,6,7 \ge \Delta$ and $d=9\ge 2\Delta -1.$
Thus the ladder is regular, see \cite{fu} Theorem 3.5. Let $\sls$ 
be  the surface section and let $\clc$ be the curve section. The
projective
normality of linearly normal surfaces of degree nine was studied in
\cite{gisa}. Comparing the list given there with \cite{fa-li9}  and
using Lemma
\ref{blowupF1} $\sls$ is seen to be projectively normal. Remark
\ref{liftingpn}
 then gives $r \xel= r \sls
= r \clc.$ 
Let now $g=5.$ Then $h^1(tL_C)=0$ for all $t\ge 1$ and from the
structural
sequence of $C$ in $\Pin{4}$ it is easy to see that $ r \clc = 3$ if and
only if
$\clc$ is $2$-normal. On the other hand \cite{fa-li9} shows that in this
case
$h^1(\calo_S)=0$ and since $h^1(L_C)=0$ it must be $h^1(L_S)=
0$ and thus $ h^2(\iof{S}{1}) = 0.$ Now the $2$-normality of
$\sls$ implies the
$2$-normality of $C$ as can be seen from $0 \to \iof{S}{1} \to
\iof{S}{2} \to
\iof{C}{2} \to 0$ and therefore $r \xel = 3.$

Let now $g=6.$ First notice that since $h^0(L_C)=5$ it is  $h^1(L_C)=1$
and thus
$0 \to \iof{C}{1} \to \oofp{4}{1} \to L_C \to 0$ shows that
$h^2(\iof{C}{1})=h^1(L_C) = 1$ i.e. $\clc$ cannot be $3$-regular.
Consider the
sequence 
\begin{equation}
\label{powersofLS}
0 \to tL_S \to (t+1)\restrict{L}{S} \to (t+1) \restrict{L}{C} \to 0
\end{equation}
for all $t \ge 1.$ Because $\deg{(t+1)\restrict{L}{C}} = 9(t+1) > 2g-2$
it is
$h^1((t+1)\restrict{L}{C}) =0$ for all $t \ge 1.$ Therefore the above
sequence gives
$h^2(t\restrict{L}{S}) =h^2((t+1)\restrict{L}{S})=0$ for all
$t\ge 1$ and thus $h^2(t\restrict{L}{S}) =0$ for all $t \ge 1$ by
Serre's Theorem.
From
\cite{fa-li9} we know that
$q(S)=0$ and
$p_g(S) = 1.$ Thus the sequence $0 \to \calo_S \to \restrict{L}{S}
\to\restrict{L}{C}\to 0$ gives
$h^1(\restrict{L}{S}) = h^2(\restrict{L}{S}) = 0.$ Then the sequence
\brref{powersofLS} for
$t=1$ gives
$h^1(2\restrict{L}{S}) =0.$ The sequence $0 \to \iof{S}{2} \to
\oofp{5}{2} \to 2\restrict{L}{S}
\to 0$  gives
$h^2(\iof{s}{2})= h^1(2L_S) =0.$ Then the sequence $0 \to \iof{S}{2} \to
\iof{S}{3}\to
\iof{C}{3} \to 0$, recalling that $\sls$ is projectively normal, gives 
$h^1(\iof{C}{3})=h^2(\iof{S}{2}) =0$, i.e. $\clc$ is $3$-normal. The
structure sequence
for $C$ in $\Pin{4}$ then easily shows that $\clc$ is $4$-regular and
thus $r \xel =
4.$

Let now $g=7.$ Noticing that $h^1(\restrict{L}{C})=2$ and recalling from
\cite{fa-li9}
that in this case $q(S) =0$ and $p_g(S) = 2$, the same argument as above
shows that
$\clc$ is not $3$-regular but it is $4$-regular thus $r \xel = 4.$

Let now $d=10.$ From \cite{fa-li10} we see that $h^1(\calo_X)=0$ and
therefore the
ladder of these manifolds is regular. Following the list given in
\cite{fa-li10} let $X$ be
a  sectional genus $6$, codimension $4$  Mukai manifold of  dimension
$3$ or $4.$
 The curve section $\clc$ is then a canonical curve in $\Pin{5}$ and as
such it is
projectively normal. Because $h^1(K_C) =1,$ $h^1(2K_C) = 0$ and the
ladder is regular, it follows from Lemma
\ref{rofpncurves} and Remark
\ref{liftingpn} that 
$r \xel = r \clc=4.$ 

 Let now $X$ be any of the remaining threefolds of degree $10$ in
$\Pin{7},$ all of
which have
$g=5,$ according to \cite{fa-li10}. Let
$\clc$ be a generic curve section. From the classification of manifolds
with
hyperelliptic section (see
\cite{BESO}) it follows that either $X$ is a hyperquadric fibration over
$\Pin{1}$ or
$C$ is not hyperelliptic. In the latter case it is $cl(C)\ge 1$ and
therefore Theorem
\ref{glcliff} gives the projective normality of $C.$ 
 Because $g=5$ it is $h^1(\restrict{L}{C})=0$ and then
Lemma \ref{rofpncurves}, the regularity of  the ladder and Remark
\ref {liftingpn} give $r \xel = r \clc= 3.$  

Let $\xel
\stackrel{\pi}{\to}\Pin{1}$ now be a  hyperquadric fibration. Consider 
$ W=
\Proj{\oofp{1}{1,1,1,1}}$ and let $\cal{T}=\calo_W(1).$  From
\cite{fa-li10} it follows that $X \in|2\cal{T} +   \pi^*(\oofp{1}{2})|$
and $L=
\restrict{\cal{T}}{X}.$ The 
higher vanishings $h^i(\iof{X}{k-i})=0$ for $i\ge 2$ required  for  the
$k$-regularity of 
$X$ are easily obtained for all $k \ge 3$ from the sequences
\begin{gather}
0\to \iof{X}{k-i}\to \oofp{7}{k-i} \to  \oof{X}{k-i} \to 0\notag \\
0\to (k-2-i)\cal{T} +\pi^*(\oofp{1}{-2}) \to (k-i)\cal{T}  \to
\oof{X}{k-i}  \to 0\notag
\end{gather}
recalling Remark \ref{leray}.

Notice that $|\cal{T}|$
embeds
$W$ in $\Pin{7}$ and the embedding is projectively normal, i.e. 
$H^0(\oofp{7}{k})
\to H^0(W, \oof{W}{k})$ is onto for all $k\ge 1.$ 
Therefore   $X$  is $k$-normal
in the embedding given by $\restrict{\cal{T}}{X},$ for some $k,$ if and
only if
$H^0(W,\oof{W}{k}) \to H^0(X,\oof{X}{k})$ is surjective and this
happens  if
and only if $H^1(W, (k-2)\cal{T} +\pi^*(\oofp{1}{-2})=0.$   

It is $H^1(W, (k-2)\cal{T}
+\pi^*(\oofp{1}{-2})) = H^1(\Pin{1},
\oofp{1}{-2} \tensor S^{k-2}\oofp{1}{1,1,1,1}).$ \\ Combining Lemma
\ref{buongg}
and the fact that 
$\mu^-(\oofp{1}{-2}
\tensor S^{k-2}\oofp{1}{1,1,1,1}) = k-4$ it is $H^1(\Pin{1},
\oofp{1}{-2} \tensor S^{k-2}\oofp{1}{1,1,1,1})=0$ for all $k \ge 3.$  On
the other
hand $H^1(\Pin{1},\oofp{1}{-2}) =1$ so $r \xel =3.$
\end{pf}

\begin{cor}
Let $\xel$ be a $n$-dimensional  polarized pair with $X$ smooth, $L$
very ample and
$\Delta
\le 3.$ Then
\begin{itemize}
\item[i)] $\xel$ is extremal if and only if it is either a hypersurface
or $\Delta=0,1.$
\item[ii)] If $\Delta=2$ then $r \xel = 3.$
\item[iii)] If $\Delta=3$ then $r \xel = 3$ unless $\xel$ is a complete
intersection of
type $(2,3)$ or a curve of genus $3$ embedded in $\Pin{3}$ as a curve of
type
$(2,4)$ on a smooth quadric hypersurface. In both these cases $r \xel =
4. $
\end{itemize}
\end{cor}
\begin{pf}
From the proof of Theorem \ref{Delta+2thm}  there are only  curves and
surfaces
with  $\Delta=3$ to consider.
If $X$ is a curve, since $c\ge 2$, it must be $g\ge 3$ and $d
\ge 6$. If
$d\ge 7$ Theorem \ref{alafujita} gives $r=3.$ If $d=6$ then $X \subset
\Pin{3}$ and 
\cite{Io1} section 7 gives three possible types for $X.$ $X$ is linked
to a twisted
cubic by two cubic hypersurfaces, $X$ is of type $(2,4)$ on a smooth
quadric or $X$
is a complete intersection of type $(2,3).$  In the first case
$\iofo{X}$ has a resolution
as in 
\brref{bordigares} for
$N=3$
  and therefore $r=3.$
In the second  case $X$ is not $2$-normal and therefore $r\ge 4.$ By
Theorem
\ref{lazran} $r \le 5.$ By \cite{GLP} $X$ cannot be extremal, therefore
$r=4.$  From
\cite{BEL2} it follows that  complete intersections of type
$(2,3)$ have $r=4.$

Assume $n=2.$  As above complete intersections of type $(2,3)$ have
$r=4.$ Following
\cite{Io1} Theorem 4.8 and section 7, it follows from Theorem
\ref{alafujita} and
Proposition
\ref{elscr} that  the only varieties  left to investigate are Bordiga
surfaces in
$\Pin{4}.$  They  have resolutions as
in \brref{bordigares} with $N=4.$ It is
straightforward to check $r=3.$
\end{pf}
\begin{cor}
Let $\xel$ be a $n$-dimensional  polarized pair with $n \ge 3,$  $X$
smooth, $L$ very
ample, 
$\Delta=4$ and $\xel$ not a hypersurface. Then $r \xel =3$ unless $\xel$
is 
a complete intersection of type $(2,2,2)$ or any threefold  in $\Pin{5}$
of degree $7$
in which cases
$r\xel =4.$
\end{cor}
\begin{pf}
Immediate from the proof of Theorem \ref{Delta+2thm}.
\end{pf}
\begin{cor}
Let $\xel$ be a $n$-dimensional  polarized pair with $n \ge 3,$  $X$
smooth, $L$ very
ample, 
$\Delta=5$ and $\xel$ not a hypersurface. Then $r \xel =3$ unless $\xel$
is 
in the following list, in
which cases
$r\xel =4.$
\begin{itemize}
\item[i)] $\xel \subset \Pin{5}$ is a threefold of degree $8$ fibered
over $\Pin{1}$ with
generic fibres complete intersections of type $(2,2)$  (see
\cite{BSS3}).
\item[ii)] $\xel \subset \Pin{6}$ is a threefold of degree $9$, $g=6$,
obtained by
blowing up a point on a Fano manifold in $\Pin{7}$, (see \cite{fa-li9}).
\item[iii)] $\xel \subset \Pin{6}$ is a threefold of degree $9$, $g=7$,
obtained by
a cubic section of a cone over the Segre embedding of $\Pin{1} \times
\Pin{2}
\subset \Pin{5}$, (see
\cite{fa-li9}).
\item[iv)] $\xel \subset \Pin{n+4}$ is a Mukai manifold of degree $10,$
$n=3,4$, $g=6$,
 (see \cite{fa-li10}).
\end{itemize}
\end{cor}
\begin{pf}
Immediate from the proof of Theorem \ref{Delta+2thm}.
\end{pf}
%
%
\section{Scrolls over curves of genus two}
\label{g2scrolls}
Let  $\xt= (\Proj{E}, \taut{E})$  be  an $n$-dimensional  scroll
over a curve $C$ of genus
$2.$ 
 From
\cite{gisa} (Lemma 5.2) it follows that $\Delta\xt = 2n$ and
$h^1(\cal{T})=0.$ These
facts will  be used without further mention throughout  this section.
The same
conclusions can also be drawn from Lemma \ref{buongg} and  the following
lemma:
\begin{lemma}
\label{muscrollg2}
Let $E$ be a rank $r$  very ample vector bundle  over a genus $2$ curve.
Then
$\mu^-(E) > 3$ and $h^1(C, S^t(E))=0$ for $t\ge 1.$
\end{lemma}
\begin{pf}
By induction on $r.$ If   $r=1$ then $E$ is semistable and very ample,
therefore
$\mu^-(E) = \mu(E) \ge 5.$
Let now $r\ge 2$ and assume $\mu^-(E) >3$ for every very ample vector
bundle of
rank up to
$r-1.$ From \cite{Io-To} it is $c_1(E) \ge 3r+1.$ If $E$ is semistable
then $\mu^-(E) =
\mu(E) =\frac{d}{r}\ge 3 + \frac{1}{r}>3.$

Let now $E$ be non semistable. Then there is a
quotient bundle $E \to Q\to 0$  such that $rk(Q) < rk(E)$ and $\mu(Q) =
\mu^-(E).$
Being a quotient of a very ample bundle on a curve, $Q$ is also very
ample and by
induction $\mu^-(E) = \mu(Q) \ge \mu^-(Q) >3.$

Because $\mu^-(S^t(E))=t \mu^-(E) > 3t \ge 3$ it is $h^1(S^t(E))=0$ from
Lemma
\ref{buongg}.
\end{pf}

\begin{prop}
\label{pnscrg2surf} Let $\xt$ be a surface scroll over a curve of genus
$g=2$ with degree
$T^2=d.$ Then
$\xt$ is projectively normal unless $d=8.$ In this case $X$ is as in
\cite{Io2} (4.2).
\end{prop}
\begin{pf}
The projective normality of such scrolls  up to degree $8$ was studied
in
\cite{alibaba} where the non projectively normal surfaces in the
statement can be
found. Let us assume $d \ge 9.$ If $E$ is semistable then $\mu^-(E) =
\mu(E) = d/2 >4$
and therefore $\xt$ is p.n. by Lemma \ref{criteriodelbutler}.
Let now $E$ be non semistable. Then $E$ admits a Harder Narasimhan
filtration  of
the form $0\to D\to E$ where $D$ is a line bundle.  Let now
$Q$ be the quotient line bundle $0\to D \to  E\to Q \to 0.$ From the
definition of
$\mu^-$ it is
$\mu^-(E)=\mu(Q) = \deg Q.$ Because $E$ is very ample so must be $Q.$ A
line
bundle $Q$  on a curve of genus $2$ is very ample if and only if $\deg Q
\ge 5.$  Thus
$\mu^-(E) >4$ and $\xt$ is p.n. by Lemma \ref{criteriodelbutler}.
\end{pf}

The above non projectively normal surface scrolls are  such because they
are not
$2$-normal, (see \cite{alibaba}). Indeed the next Proposition and Lemma
\ref{muscrollg2} show that
$2$-normality is equivalent to projective normality for scrolls of genus
$2,$ 
extending a result found in
\cite{pu-ga}.
\begin{prop}
\label{pn2n}
Let $\xt=(\Proj{E},\taut{E})$ be a scroll over a smooth curve $C$ of
genus $g$ such
that
$\mu^-(E)>2g-2.$
Then $\xt$ is projectively normal if and only if it is $2$-normal.
\end{prop}
\begin{pf}
If $\xt$ is p.n. it is obviously $2$-normal. \\ To see the converse
let $n=dim X=rk E$ and let $\pi:X\to C$ be the natural projection.
Reasoning as in \cite{pu-ga} Lemma 1.4, projective normality of
$\xt$ follows from  the surjectivity of the maps 
\begin{equation}
\label{ontomaps}
H^0((j-1)\cal{T}) \tensor H^0(\cal{T}) \to
H^0(j\cal{T}) \text{ \ \ \ for all \ }
j \ge 2.
\end{equation}  
This in turns follows, according to \cite{mu} Theorem 2, from the
vanishing of
$$H^i(X,(j-1-i)\cal{T})=0 \text{\ \    for all } n\ge i\ge 1 \text{ \
\      and for all } j \ge
2.$$ Because $i\le n$ and $j\ge 2$ it is $j-1-i>-n$ and therefore
Remark \ref{leray} shows  that $\xt$ is p.n. if
$H^1(X,(j-2)\cal{T})=H^1(C,
S^{j-2}E)=0$ for all $j\ge 2.$
The hypothesis $\mu^-(E)>2g-2$ implies $\mu^-(S^{j-2}E) =
(j-2)\mu^-(E)>2g-2$
for all $j\ge 3.$ From  Lemma \ref{buongg} it follows that
$H^1(X,(j-2)\cal{T})=0$ for
all
$j\ge 3.$ This gives all necessary surjectivity in \brref{ontomaps} but
for $j=2.$
Thus $\xt$ is p.n. if
$H^0(\cal{T}) \tensor H^0(\cal{T}) \to H^0(2\cal{T})$ is onto, i.e. if
$\xt$ is $2$-normal. 
\end{pf}
\begin{cor}
A scroll $\xt$ over a curve of genus $2$ is p.n. if and only if  it is
$2$ normal.
\end{cor}
\begin{pf}
Immediate from Proposition \ref{pn2n} and Lemma \ref{muscrollg2}
\end{pf}
Results on threefold scrolls are collected in the following proposition.
\begin{prop}
Let $\xt=(\Proj{E}, \taut{E})$ be a threefold scroll over a curve of
genus $2.$ Let $d
\ne 12.$ Then
$\xt$ fails to be projectively normal if and only if one of the
following cases occur
\begin{enumerate}
\item[i)] $d=11.$
\item[ii)] $d\ge 13$, $E$ is not semistable and it admits  a quotient $E
\to \cal{E}\to 0$ of
rank two and degree eight.
\end{enumerate}
\end{prop}
\begin{pf}
It is known, see \cite{Io-To} or \cite{Io2} and \cite{fa-li9},
\cite{fa-li10},
that there do not
exist threefold scrolls of genus two and
$d\le 10.$ So assume $d\ge 11.$ Because $h^1(\cal{T})=0$ it is
$h^0(\cal{T}) = d-3.$ A simple
computation shows that
$h^0(X,\oof{X}{2}) = 4d-6 > h^0(\oofp{d-4}{2}) = \frac{(d-2)(d-3)}{2}$
if $d \le 11,$ so
that degree $11$ scrolls cannot be $2$-normal.
Assume $d\ge 13.$ If $E$ is semistable then $\mu^-(E) = \mu(E) = d/3 >
4$ and thus
$\xt$ is p.n. by Lemma \ref{criteriodelbutler}.
Let now $E$ be not semistable. Assume $E$ does not admit  a degree $8$
and rank $2$
quotient. All  quotients $E \to Q \to 0$ must be very ample and thus it 
must be $rank
Q = 1$ and
$\deg Q \ge 5$ or $rank Q =2$ and $\deg Q \ge 9.$ Therefore for all $Q$
it is  $\mu(Q)
>4$ and thus $\mu^-(E) >4$ and then $\xt$ is p.n. by Lemma \ref{criteriodelbutler}.

Let now $E$ be not semistable with a quotient $E\to \cal{E} \to 0$ of
degree $8$ and
rank $2.$ Notice that $(\Projcal{E},\tautcal{E})$ is one of the non
$2$-normal
surfaces of degree eight embedded in $\Pin{5}$ studied in
\cite{alibaba}.
Let $D$ be the line bundle of degree $d -8$  such that 
\begin{equation}
\label{quotient}
0\to D \to  E \to
\cal{E}
\to 0.
\end{equation}
Since $d - 8 > 2$ it  is $h^1(D) =0$ and thus $H^0(E) = H^0(\cal{E})
\oplus H^0(D).$
Therefore 
\begin{equation} 
\label{essedueE}
S^2(H^0(E))=S^2(H^0(\cal{E})) \bigoplus H^0(\cal{E}) \tensor
H^0(D)\bigoplus
S^2(H^0(D))
\end{equation}
Consider the sequence obtained by tensoring \brref{quotient} with $D:$
\begin{equation} 
\label{tensorD}
0 \to D \tensor D \to E \tensor D \to \cal{E} \tensor D \to 0.
\end{equation}
Because $\deg(D\tensor D) =2(d-8) > 2$ and $\mu^-(\cal{E} \tensor
D)=\mu^-(\cal{E})
+ \mu^-(D) = d-4 > 2$ it follows that $h^1(D\tensor D) = h^1( \cal{E}
\tensor D)=0$ and
thus $H^0(E \tensor D) = H^0(D\tensor D) \oplus H^0(\cal{E} \tensor D)$
and
$h^1(E \tensor D) = 0.$
Considering  now the exact sequence
\begin{equation}
\label{simm}
0 \to D \tensor E \to S^2(E) \to S^2(\cal{E}) \to 0
\end{equation}
it follows that 
\begin{equation}
\label{hzero}
H^0(S^2(E)) = 
H^0(S^2(\cal{E})) \bigoplus H^0(\cal{E} \tensor D)\bigoplus H^0(D
\tensor D).
\end{equation}
Putting together \brref{hzero} and \brref{essedueE}
 it follows that the map $\phi: S^2(H^0(E)) \to H^0(S^2(E))$ decomposes
as 
\begin{gather}[S^2(H^0(\cal{E}))
\stackrel{\alpha}{\rightarrow}H^0(S^2(\cal{E}))] \notag \\ \bigoplus 
\notag \\
[H^0(\cal{E})
\tensor H^0(D) \stackrel{\beta}{\rightarrow} H^0(\cal{E}\tensor D)]
\bigoplus [S^2(H^0(D)) \stackrel{\gamma}{\rightarrow} H^0(D \tensor D)].
\notag
\end{gather}

It was proven in \cite{alibaba} that $\alpha$ is not surjective,
therefore $\phi$
cannot be surjective, i.e. $\xt$ cannot be $2$-normal.
\end{pf}
\begin{rem}
The existence of degree $11$ and $12$ threefold scrolls over curves of
genus $2$ 
is an open problem. If a degree $12$ such scroll
$\xt=(\Proj{E},\taut{E})$ exists then it is not difficult to see that
$E$ must be
semistable. If it were not semistable then there would be  a
destabilizing  subbundle\ 
$\cal{F}$ either of rank $1$ such that
$ \deg \cal{F} > 4$ or of rank $2$ such that $ \deg
\cal{F} > 8.$ 

In both cases the resulting quotient $0 \to \cal{F} \to E \to Q \to 0$
could
not be very ample for degree reasons, which is a contradiction.
\end{rem}
\begin{rem} Let $\xt=(\Proj{E},\taut{E})$ again be a $3$ dimensional
scroll over a
curve of genus $2.$ If $\xel$ is projectively normal, recalling that
$h^1(E) =0,$
the same argument used in Proposition \ref{elscr} gives $r \xel = 3.$ 
If  $\xel$ is not p.n. , notice that if
$d\ge 13$ it follows that $h^0(\cal{T}) \ge 10$ and thus $(X, \cal{T})$
is $(\Delta +
2)$-regular, i.e. $8$-regular  from  \cite{Ran1}.
When $d=11,12$ it is easy to check that $h^i(\iof{X}{8-i}) = 0$ for all
$i\ge 2$ while
we were not able to  establish the $7$-normality of these manifolds. 
\end{rem}
%
%
\section{ $\Pin{r}$ bundles over an elliptic curve }
\label{ellipticpkbundles}
Throughout this section let $E$ be a vector bundle of rank $r$ and
degree $d$ over an
elliptic curve
$C$. Let
$\xt =
\scroll{E}$ and let $D$ be a divisor on $X$ numerically equivalent to
$aT + bF.$ 
Assume $D$ is very ample. The projective normality of the embedding
given
by $D$ was studied by Homma
\cite{Ho1},
\cite{Ho2} when $r=2,$  and in a more general setting by Butler
\cite{bu} (See also
\cite{alibaba}). In this section the case of $a=2$ and $r=3$ is
addressed and  Butler's results are improved in some cases.

\begin{lemma}
\label{mupiudec}
Let $E = \bigoplus_i E_i$ be a decomposable vector bundle over an
elliptic curve. Then
$$\mu^+(E) = max_i \{{\mu(E_i)}\}.$$
\end{lemma}
\begin{pf}
This is essentially   \cite{alibaba2} Lemma 2.8, reinterpreted from the
point of view
of $\mu^+$ instead of $\mu^-.$
\end {pf}

\begin{lemma}
\label{mupiu}
Let $\xt$ be as above and let $M_{s}$ be a divisor on $X$ whose
numerical class is
$M_{s}\equiv T+sF. $ Let
$m=min\{ t 
\in \Bbb{Z}| h^0(M_t)> 0.\}$ Then $m=-[\mu^+(E)]$ and there exists a
smooth $S \num T+mF.$
\end{lemma}
\begin{pf}
From \cite{bu} it follows that for any vector bundle $\cal{G}$ over an
elliptic
curve $\mu^+(\cal{G}) <0$ implies $h^0(\cal{G}) = 0.$
For simplicity of notation let  $m^*=-[\mu^+(E)].$ We need to show that
$m=m^*.$  Let
$\cal{L}_t$ be a line bundle on $C$ with degree $t.$ If $t < m^*,$ then
$t=m^*- x$ for some integer
$x\ge 1.$ Then
$\mu^+(E \tensor\cal{L}_t )= \mu^+(E) + t = \mu^+(E) +m^*- x < 1-x \le
0$
Therefore
$ h^0(M_t) = 0$ if $t< m^*$ and thus 
\begin{equation} \label{mandm1} m^*\le m.\end{equation}

Let now $E$ be indecomposable and thus semistable. Because $\mu(E)
=\mu^+(E)$ and
$-m^*
\le
\mu^+(E)$ it is $d+rm^* \ge 0.$  If $d +rm^*>0$ then $
h^0(M_{m^*}) >0.$ If $d + rm^*= 0$ then, as in \cite{At}, a  line
bundle
$\cal{L}_{m^*}$ of degree $m^*$ can be found by a suitable twist of
degree zero,  such
that
$h^0(E \tensor \cal{L}_{m^*}) =1.$

Let now $E = \bigoplus_i E_i$ be decomposable. Then $h^0(E) =
\bigoplus_ih^0(E_i).$
Let $E_{\hat{i}}$ be one of the components such that $\mu^+(E) =
\mu(E_{\hat{i}}).$
As $E_{\hat{i}}$ is indecomposable  it follows from above that  there
exists a  line
bundle
$\cal{L}_{m^*}$ of degree $m^*$  such that 
$h^0(E_{\hat{i}}
\tensor
\cal{L}_{m^*})>0.$ 

Therefore
$h^0(M_{m^*}) > 0$ and thus $m\le m^*$ which combined with
\brref{mandm1} gives $m=m^*.$

 If $S\num T+m^*F$ is singular
it must be reducible as
$S'
\cup (m^*-t) F$ where
$S'\num T+tF$ with $t<m^*$ which is not possible because of the
minimality of
$m^*=m.$ Therefore there must be a smooth $S\num T +m^*F.$
\end{pf}
\begin{lemma}
\label{D2npn}
Let  $\xt$  and $D$ be as above with $r=3,$ $a \ge 2,$ and $D$ very
ample. If
\begin{itemize}
\item[i)] there exists an ample  smooth surface $S \num T+x F$ for some
$x \in
\Bbb Z;$ 
\item[ii)] $(a-1)\mu^-(E) + b - x>1$
\end{itemize}
then
the embedding of $X$ given by $D$ 
 is projectively normal if and only if it is 
$2$-normal. 
\end{lemma}
\begin{pf}
If the embedding is p.n. it is obviously $2$-normal. As in Lemma
\ref{pn2n} the
projective normality follows from  the surjectivity of the maps 
\begin{equation}
\label{ontomaps2}
H^0((j-1)D) \tensor H^0(D) \to
H^0(jD) \text{ \ \ \ for all \ }
j \ge 2.
\end{equation} 
Assume $j\ge 4.$ 
Surjectivity in \brref{ontomaps2}  follows, according to \cite{mu}, from
the
vanishing of
$$H^i(X,(j-1-i)D)=0 \text{\ \    for all } 3 \ge i\ge 1 \text { \ \     
and for all } j \ge 4.$$
Notice that  $R^q\pi_*((j-1-i)D)=R^q\pi_*(a(j-1-i) \cal{T}) \tensor
\cal{M}_{i,j}$ where
$\cal{M}_{i,j}$ is a line bundle on $C$ of degree $(j-1-i)b.$  Notice
also, e.g. \cite{H}, that
$R^q\pi_*(a(j-1-i)\cal{T})=0$ unless
$q=0$ and $j-1-i\ge 0,$ or $q=2$ and $a(j-1-i)\le -3.$ Since $a\ge 2,$ 
$i\le 3$ and
$j\ge 4,$ the last inequality is never satisfied . For $j \ge 4$
Leray's spectral sequence shows that it is enough to show
$H^1(X,(j-2)D)=H^1(C, S^{a(j-2)}E\tensor \cal{M}_{1,j})=0$ which is
guaranteed by 
$D$ being ample.  This gives all necessary surjectivity in
\brref{ontomaps2} but
for
$j=2,3.$
If $j=2$ the map $H^0(D) \tensor H^0(D) \to H^0(2D)$ is onto by
assumption , being 
the embedding  $2$-normal. 
Assume now  $j=3.$ 
Let $S  \num T+x F$ be the smooth element whose existence is given by
assumption $i).$  Ampleness of $D$ gives $a \mu^-(E) + b > 0.$ Combining
this with
condition $ii)$ it follows from Lemma \ref{buongg} that $H^1(tD-S) =0$
for $t=1,2,3.$ 
In particular, following Homma,  the commutative diagram below is
obtained :

\begin{alignat}{5}
0 \to H^0(D-S)&\tensor H^0(2D)& &\to& H^0(D)&\tensor H^0(2D) &&\to&
H^0(S,\restrict{D}{S})
&\tensor H^0(2D)\to 0 \notag \\ 
&\downarrow \alpha& & &  &\downarrow \beta &   &  &  &\downarrow
\gamma\\ 
0 \longrightarrow H^0(3&D-S)& &\longrightarrow& H^0(3&D)&
&\longrightarrow& H^0(S,&\restrict{3D}{S}) \to 0  \notag
\end{alignat}
The surjectivity of $\beta$ will follow from the surjectivity of
$\alpha$ and $\gamma.$
From \cite{Ho1} and \cite{Ho2} it follows that  $\restrict{D}{S}$ is
normally generated.
Because
$H^0(2D)
\to H^0(\restrict{2D}{S})$ is surjective from above, it follows that
$\gamma$ is onto.

Lemma \ref{buongg} and condition $ii)$ give $D-S\num (a-1)T + (b-x)F$ 
generated
by global sections. Using this fact and noticing that $H^1(D+S)=0$ being
$D$ very
ample and
$S$ ample, it is straightforward to check that $H^i(2D - i (D-S)) = 0$
for all $i \ge 1.$
Therefore by
\cite{mu} $\alpha$ is onto.
\end{pf}
\begin{prop}
\label{homma}
Let $\xt$ and $D \num 2T + bF$ be as above.  If 
\begin{itemize}
\item[i)] there exists an ample  smooth divisor $Y \num T+x F$ for some
$x \in
\Bbb Z;$ 
\item[ii)] $\mu^-(E) + b - x>1$
\end{itemize}
then $|D|$ gives a $2$-normal embedding of $X$  if  $|\restrict {D}{Y}|$
gives a
$2$-normal embedding of $Y$.
\end{prop}
\begin{pf}
The proof proceeds along the same lines of the case $j=3$ in the proof
of Lemma
\ref{D2npn}. Let
$Y \num T+x F$ be the smooth element whose existence is given by
assumption
$i).$ Ampleness of $D$ gives $2 \mu^-(E) + b > 0.$ Combining this with
condition $ii)$ it
follows that $H^1(tD-Y) =0$ for $t=1,2.$  In particular the following 
commutative diagram  is obtained :

\begin{alignat}{5}
0 \to H^0(D-Y)&\tensor H^0(D)& &\to& H^0(D)&\tensor H^0(D) &&\to&
H^0(Y,\restrict{D}{Y})
&\tensor H^0(D)\to 0 \notag \\ 
&\downarrow \alpha& & &  &\downarrow \beta &   &  &  &\downarrow
\gamma\\ 
0 \longrightarrow H^0(2&D-Y)& &\longrightarrow& H^0(2&D)&
&\longrightarrow& H^0(Y,&\restrict{2D}{Y}) \to 0.  \notag
\end{alignat}
The surjectivity of $\beta$ will follow from the surjectivity of
$\alpha$ and $\gamma.$

Because $H^0(D) \to
H^0(\restrict{D}{Y})$ is onto from above and $H^0(\restrict{D}{Y})
\tensor
H^0(\restrict{D}{Y}) \to H^0(\restrict{2D}{Y})$  is onto by assumption
it follows that
$\gamma$ is onto.

Condition $ii)$ is equivalent to $D-S\num T + (b-x)F$ being generated by
global
sections (see Lemma \ref{buongg} and \cite{alibaba2} Lemma 2.9). Using
this fact and
noticing that
$H^1(Y)=0$ being
$Y$ ample, it is straightforward to check that $H^i(D - i(D-Y)) = 0$ for
all $i \ge 1.$
Therefore by
\cite{mu} $\alpha$ is onto.
\end{pf}

\begin{cor}
\label{hommacor}
Let $\xt$ and $D$ be as above with $r=3$ and $a=2.$ If 
\begin{itemize}
\item[i)]$\mu^-(E) > [\mu^+(E)]$
\item[ii)] $\mu^-(E) + b >1-[\mu^+(E)]$
\end{itemize}
then $|D|$ gives a projectively normal embedding.
\end{cor}
\begin{pf}
Let  $x=-[\mu^+(E)]$, notice that condition i) and Theorem
\ref{miyaoteo}
give ampleness of $T+xF$. Now combine Lemma
\ref{mupiu}, Lemma
\ref{D2npn},  Proposition
\ref{homma} and the fact that a very ample line bundle on a $2$
dimensional scroll
over an elliptic curve is always normally generated by \cite{Ho1} and
\cite{Ho2}.
\end{pf}
\begin{rem}
Let $E$ be an indecomposable vector bundle of rank $r=3$ and degree
$d\num 1 (3).$
For simplicity let us assume that $E$ has been normalized, so $d=1.$
Since $E$ is
indecomposable it is semistable and $\mu^-(E) = \mu ^+(E) = \mu(E) =
1/3.$ The
hypothesis of Corollary \ref{hommacor}  are satisfied for $D\num 2T +
F.$  
and such a $D$ is very ample from \cite{alibaba2} Theorem 4.5. Therefore
$|D|$ gives
a projectively normal embedding.  Notice that  Butler's results
\cite{bu}, see Lemma
\ref{criteriodelbutler}, were not able to establish the normal
generation of such a $D.$
\end{rem}
\begin{rem}
It is straightforward to check that for a divisor  $D$ as in Corollary
\ref{hommacor}  it is always $h^0(D) \ge 10$ and therefore the embedding
given by
$|D|$ satisfies the Eisenbud Goto conjecture .
\end{rem}

\end{document}